\documentclass[journal]{IEEEtran}

\usepackage{geometry}
\usepackage{fullpage}
\usepackage{amssymb}
\usepackage{amsfonts}
\usepackage{xcolor}
\usepackage{graphicx}
\usepackage{hyperref}
\hypersetup{colorlinks=true}
\usepackage{xcolor}
\usepackage{amsmath}
\usepackage{algorithmic,algorithm}
\usepackage{relsize,mathtools,nccmath}

\newcommand {\R}   {{\rm I\!R}}

\DeclareMathOperator*{\argmax}{arg\,max}
\DeclareMathOperator*{\argmin}{arg\,min}

\newcommand{\bfK}{{\bf K}}

\newcommand{\bfI}{{\bf I}}

\newcommand{\bfc}{{\bf c}}

\newcommand{\bfa}{{\bf a}}

\newcommand{\bfy}{{\bf  y}}
\newcommand{\bfx}{{\bf  x}}

\newcommand{\bfz}{{\bf z}}

\newcommand{\Ex}{\mathbb{E}}
\newcommand{\Lc}{\mathcal{L}}
\newcommand{\Jc}{\mathcal{J}}

\newcommand{\Lcurly}{\mathcal{L}}

\newcommand{\bfzeta}{\boldsymbol{\zeta}}
\newcommand{\bfmu}{\boldsymbol{\mu}}
\newcommand{\bfSigma}{\boldsymbol{\Sigma}}

\newcommand{\bftheta}{\boldsymbol \theta}

\title{Kernel Expansions for High-Dimensional Mean-Field Control with Non-local Interactions}
\author{Alexander Vidal, Samy Wu Fung, Stanley Osher, Luis Tenorio, Levon Nurbekyan}

\begin{document}

\maketitle
%

\begin{abstract}
Mean-field control (MFC) problems aim to find the optimal policy to control massive populations of interacting agents. 
These problems are crucial in areas such as economics, physics, and biology. 
We consider the non-local setting, where the interactions between agents are governed by a suitable kernel. For $N$ agents, the interaction cost has $\mathcal{O}(N^2)$ complexity, which can be prohibitively slow to evaluate and differentiate when $N$ is large.
To this end, we propose an efficient primal-dual algorithm that utilizes basis expansions of the kernels. The basis expansions reduce the cost of computing the interactions, while the primal-dual methodology decouples the agents at the expense of solving for a moderate number of dual variables. 
We also demonstrate that our approach can further be structured in a multi-resolution manner, where we estimate optimal dual variables using a moderate $N$ and solve decoupled trajectory optimization problems for large $N$. We illustrate the effectiveness of our method on an optimal control of 5000 interacting quadrotors.


\end{abstract}


\IEEEpeerreviewmaketitle

\section{Introduction}
Mean-Field Control (MFC)~\cite{bensoussan2013mean} is a class of optimal control problems that primarily focus on controlling systems comprising a large number of identical interacting agents. These problems are designed to  optimize  \textit{collective} rather than individual behaviour of the agents; in other words, the goal is to control the statistical properties of the population, e.g., the distribution of the agents in the state space or the physical space. Therefore, MFC problems can be found in numerous applications, including epidemic modeling~\cite{aurell2022optimal, lee2022mean, charpentier2020covid}, finance~\cite{cardaliaguet2018mean, djehiche2015stochastic, calvia2024mean}, and water distribution~\cite{barreiro2019mean, barreiro2021mean}. 
\begin{figure}[t]
    \centering
    \includegraphics[width=0.28\textwidth]{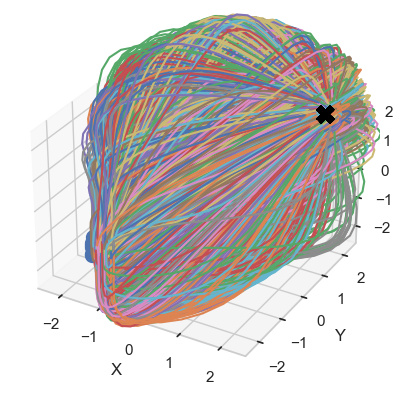}
    \caption{Optimal swarm trajectories for 5000 quadrotors using our proposed primal-dual approach. The black X denotes the target. }
    \label{fig: 5k_quadcopter_trajectories}
\end{figure}
There are two primary computational challenges when solving MFC problems, the curse-of-dimensionality, and the coupling among agents arising from \emph{non-local interactions}. 
While many methods have addressed the curse-of-dimensionality when solving MFCs~\cite{ruthottomean2020, HanSolving2018, zhao2022offline, lin2021alternating}, 
a key computational challenge remains when solving high-dimensional MFC problems with \emph{non-local} interactions; 
namely, when modeling $N$ agents with non-local interactions~\cite{Nurbekyan2019Fourier}, evaluation of the interaction cost has complexity $\mathcal{O}(N^2)$ \emph{at any given time}, rendering existing MFC approaches computationally expensive (see Fig.~\ref{fig: interaction_cost_time}).

\subsection{Our Contribution}
We propose a primal-dual framework based on kernel expansions to solve high-dimensional MFC problems with \emph{non-local interactions}. In particular, by introducing coefficients that capture the mean-field interaction~\cite{NurbekyanOne2018}, our proposed primal-dual approach decouples all the agents at the expense of the updates of a moderate number of dual variables, rendering the optimization of primal variables parallelizable. In our numerical results, we are able to approximate optimal control of 5000 quadrotor agents (see Figure~\eqref{fig: 5k_quadcopter_trajectories}). Accompanying animated graphics and code for our experiments can be found in \url{https://github.com/mines-opt-ml/kernel-expansions-for-mfc}.

\section{Background}

We are interested in MFC problems of the form
\begin{equation}
    \begin{split}
    \inf_{\bftheta, \rho} \,\,\, & \Jc(\bftheta, \rho) \\
    \text{s.t. } & \partial_t \rho_t(\bfx) + \nabla \cdot (\rho_t(\bfx) f(t,\bfx, \bftheta) ) = 0, \quad t \in (0,T), 
    \end{split}
    \label{eq: original_MFC_problem}
\end{equation}
where $\{\rho_t\}$, $t\in [0,T]$, are probability densities on $\mathbb{R}^d$ with known
initial density $\rho_0$, representing the distribution of the agents in the state-space at each time $t$. The objective function
$\mathcal{J}$ is defined as
\begin{equation}
    \begin{split}
        \Jc(\bftheta, \rho) = & \medmath{\int_0^T}\left[\, \Ex_{\bfx\sim\rho_s}L(s,\bfx,\bftheta_\bfx(s))  + \mathcal{F}(\rho_s) \,\right] ds\\
        + & \,\,\mathcal{G}(\rho_T),
        \end{split}
        \label{eq: mfc_original_objective}
\end{equation}
where $\bftheta:\R^d\times [0,T] \to \mathcal{U} \subset \R^q$ is the control and $\mathcal{U}$ is a sufficiently regular domain (see~\cite{fleming2006controlled} see Sec. I.3, I.8-9);  
the function $f : [0, T] \times \R^d \times \R^q \to \R^d$ models the dynamics of the agents, and $L \colon [0,T] \times \mathbb{R}^d \times \mathcal{U} \to \mathbb{R}$ is the Lagrangian that controls a particular performance criterion, e.g., kinetic energy. Note that the indistinguishability of the agents yields controls that depend (besides the time variable) only on the position of one agent at a time.

The functionals $\mathcal{F}$ and $\mathcal{G}$ model the interactions among agents within the time-horizon $(0,T)$ and at the final time $T$, respectively. For the terminal interaction cost we use
\[
\mathcal{G}(\rho_T) = \Ex_{\bfx\sim\rho_T} G(\bfx)
\]
for a suitably chosen $G$.

For the non-local interactions within the time horizon $(0,T)$ we use
\begin{equation}
    \mathcal{F}(\rho_t) = \frac{1}{2}\Ex_{(\bfx,\bfy)\sim \rho_t\otimes\rho_t} K(\bfx,\bfy),
\label{eq: mean_field_interaction}
\end{equation}
where $K\colon \mathbb{R}^d \times \mathbb{R}^d \to \mathbb{R}$ is a positive definite kernel~\cite[Chapter 6]{mohri18foundations}.
In our examples, the interaction term penalizes agents from being too close to one another.
Thus, we want $K(\bfx,\bfy)$ to be large when $\bfx$ and $\bfy$ are close and small otherwise. Common choices of kernels include a Gaussian kernel~\cite{Nurbekyan2019Fourier, Agrawal2022RandomFF} or the inverse distance kernel~\cite{lin2021alternating, fasshauer2015kernel}.
\begin{figure}[H]
    \centering
    \includegraphics[width=0.3\textwidth]{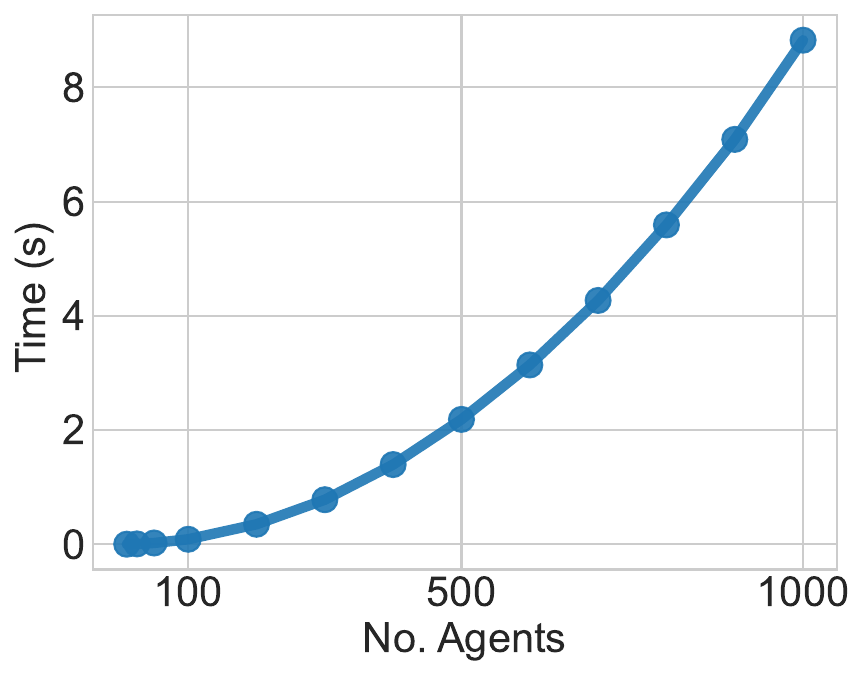}
    \caption{\textbf{Non-local mean-field interaction cost per evaluation:} Average evaluation timings for mean-field interaction term~\eqref{eq: mean_field_interaction}.  Averages are obtained by evaluating the mean-field interaction term three times and taking the sample average of the evaluation times.}
    \label{fig: interaction_cost_time}
\end{figure}
\noindent 
The double expectation (or integral) in~\eqref{eq: mean_field_interaction} results in a significant computational burden when solving the MFC problem~\eqref{eq: mfc_original_objective} because evaluation of $\mathcal{F}$ for $N$ agents results in $\mathcal{O}(N^2)$ complexity, see Figure~\ref{fig: interaction_cost_time}.

\section{{Primal-Dual Framework for Agent Decoupling}}
\subsection{Saddle Point Formulation}
\label{sec: PD via Kernel Expansions}
To address the computational challenges associated with evaluation of the interaction cost~\eqref{eq: mean_field_interaction}, we extend the kernel expansion framework introduced in \cite{Liucomputational2021, Liusplitting2021, NurbekyanOne2018, Nurbekyan2019Fourier, rahimiadvances2007} to  MFC problems.  
We begin by approximating the kernel in~\eqref{eq: mean_field_interaction} with with a quadratic (in the feature space) function
\begin{equation}\label{eq: kernel_expansion_framework}
K(\bfx,\bfy) \approx  K_r(\bfx,\bfy) = \bfzeta(\bfx)^\top\bfK_r\,\bfzeta(\bfy),
\end{equation}
where $\bfzeta = (\zeta_1, \zeta_2, \ldots, \zeta_r)^\top$,
$\zeta_i\in \mathcal{C}^2(\R^d)$ are some basis functions or features, and $\bfK_r$ is an
$r\times r$ symmetric positive-definite matrix. Note that when $\bfK_r$
is diagonal we obtain a truncated Mercer series approximation \cite{mercer1909xvi}. 
 While there are many ways to choose $\bfzeta$~\cite{rahimiadvances2007, Agrawal2022RandomFF}, in this work we use a neural network approximation. More details are provided in Sec.~\ref{sec: numerical_experiments}. The objective function \eqref{eq: mfc_original_objective} can now be approximated using
\begin{equation}
    \begin{split}
      &\Jc_r(\bftheta, \rho) = \medmath{\int_0^T} \Ex_{\bfx\sim\rho_s}L(s,\bfx,\bftheta_\bfx(s)) \,ds \\
      &+ \frac{1}{2}\medmath{\int_0^T} \Ex_{(\bfx,\bfy)\sim\rho_s\otimes\rho_s} K_r(\bfx,\bfy) \,ds\\
      &+ \mathcal{G}(\rho_T).
      \label{eq: mean_field_control_problem_Kr}
    \end{split}
\end{equation}

Next, we rewrite the approximate interaction term in~\eqref{eq: mean_field_control_problem_Kr} using \eqref{eq: kernel_expansion_framework} to obtain 
\begin{equation}
\label{eq: kernel_approximation}
\begin{split}
      & \frac{1}{2}\Ex_{(\bfx,\bfy)\sim\rho_t\otimes\rho_t} K_r(\bfx,\bfy)
    \\
    &=\frac{1}{2}\Ex_{\bfx\sim\rho_t}\bfzeta(\bfx)^\top\,
    \mathbf{K}_r \Ex_{\bfy\sim\rho_t} \bfzeta(\bfy)  
      \\
      &=  \frac{1}{2}\bfc_t^\top \mathbf{K}_r \bfc_t,
\end{split}
\end{equation}
where
\begin{equation*}
    \bfc_t = \Ex_{\bfx\sim\rho_t} \bfzeta(\bfx),\quad t\in (0,T).
\end{equation*}
Using convex duality, we have that
\begin{equation}\label{eq:duality}
\begin{split}
    \frac{1}{2}\bfc_t^\top \mathbf{K}_r \bfc_t=&\sup_{\bfa_t} \left[ \bfa_t^\top \bfc_t-\frac{1}{2} \bfa_t^\top \mathbf{K}_r^{-1} \bfa_t \right]\\
    =&\sup_{\bfa_t} \left[ \Ex_{\bfx\sim\rho_t} \bfa_t^\top \bfzeta(\bfx)-\frac{1}{2} \bfa_t^\top \mathbf{K}_r^{-1} \bfa_t \right].
\end{split}
\end{equation}
Hence, optimizing~\eqref{eq: mean_field_control_problem_Kr} is equivalent to the saddle point problem
\begin{align*}
        &\inf_{\bftheta,\rho} \sup_\bfa \medmath{\int_0^T} \Ex_{{\bfx\sim\rho_s}} L(s,\bfx,\bftheta_\bfx(s))\,ds + \mathcal{G}(\rho_T) \notag  \\
        &\,-\, \medmath{\int_0^T} \Ex_{\bfx\sim\rho_s}\bfa_s^\top \bfzeta(\bfx)\, ds -\medmath{\frac{1}{2} \int_0^T }\bfa_s^\top\, \bfK_r^{-1}\,\bfa_s\, ds.
\end{align*}
Interchanging the $\inf$ and $\sup$ and changing the sign of the objective function, we arrive at the variational formulation in~\cite{Nurbekyan2019Fourier}:
\begin{align}
        &\inf_\bfa \sup_{\bftheta}\medmath{\frac{1}{2} \int_0^T }\bfa_s^\top\, \bfK_r^{-1}\,\bfa_s\, ds \notag \\
        &- \medmath{\int_0^T} \Ex_{{\bfx\sim\rho_s}} L(s,\bfx,\bftheta_\bfx(s))\,ds  \notag
        - \medmath{\int_0^T} \Ex_{\bfx\sim\rho_s}\bfa_s^\top \bfzeta(\bfx)\, ds \notag \\
        &\,-\, \mathcal{G}(\rho_T) \notag \\
        &\text{s.t. } \partial_s \rho_s(\bfx) + \nabla \cdot (\rho_s(\bfx) f(s,\bfx,\bftheta_\bfx(s))) = 0, ~s \in (0,T). \notag\\
    \label{eq: saddle_point_formulation}
\end{align}
Since Lagrangian coordinates are more suitable for high-dimensional problems~\cite{ruthottomean2020, lin2021alternating, HanDeepLA2016, HanSolving2018}, we reformulate~\eqref{eq: saddle_point_formulation} in Lagrangian coordinates. More specifically, let
\begin{equation}
    \partial_s \bfz_\bfx(s) = f(s, \bfz_\bfx(s), \bftheta_\bfx(s)), \quad \bfz_\bfx(0) = \bfx,
    \label{eq : dynamics}
\end{equation}
for $s \in (0,T)$, where we suppress the dependence of $\bfz$ on $\bftheta$ until it is necessary.

Then $\bfx \sim \rho_0$ yields $\bfz_\bfx(s)\sim \rho_s$, where $\rho_s$ is the solution of the constraints equation in~\eqref{eq: saddle_point_formulation}. Consequently,~\eqref{eq: saddle_point_formulation} can be written as
\begin{equation}
    \begin{split}
    &\inf_\bfa \sup_{\bftheta} \Ex_{\bfx\sim \rho_0} \; \Big[\medmath{\frac{1}{2} \int_0^T} \bfa_s^\top \,\bfK_r^{-1}\,\bfa_s\,ds \,\,\,\,\,\\
    &-\medmath{\int_0^T}  L(s,\bfz_\bfx(s),\bftheta_\bfx(s))\,ds - \medmath{\int_0^T} \bfa_s^\top\bfzeta(\bfz_\bfx(s))\, ds \\
    & - G(\bfz_\bfx(T))\Big],\\
    \text{s.t.~}&\eqref{eq : dynamics}~\text{holds.}
    \end{split}
    \label{eq: lagrangian_coordinates_saddlepoint}
\end{equation}
The function $\bfa_t$ encodes information about the interactions among the agents. Indeed, the optimization over $\bftheta$ in~\eqref{eq: lagrangian_coordinates_saddlepoint} \textit{decouples} into regular optimal control problems
\begin{equation}
\begin{split}
        \phi_{\bfa}(t,\bfx) =& \inf_{\bftheta} \medmath{\int_t^T } L(s,\bfz_x(s),\bftheta_\bfx(s)) ds\\
        &+  \medmath{\int_t^T} \bfa_s^\top \bfzeta(\bfz_\bfx(s))\, ds + G(\bfz_\bfx(T)),\\
        \text{s.t.~}&\eqref{eq : dynamics}~\text{holds.}
\end{split}
\label{eq: phi_a_def}
\end{equation}
Thus, once an \textit{optimal} $\bfa_t$ is known, agents only need to solve decoupled optimal control problems, and the interaction information is reflected in the dependence of the cost function in~\eqref{eq: phi_a_def} on $\bfa_t$. Hence, throughout this paper, we refer to $\bfa_t$ as the \emph{global interaction coefficients.}

The function $\phi_\bfa$ is called the value function~\cite{fleming2006controlled}; it is the viscosity solution of the Hamilton-Jacobi equation 
\begin{equation}
    \begin{cases}
    -\partial_t \phi(t,\bfx) + H(t, x, \nabla_\bfx \phi(t,\bfx)) = \bfa_t^\top\bfzeta(\bfx)\\ 
    \phi(T,\bfx) = G(\bfx).
\end{cases} 
\end{equation}
Eliminating the optimization over $\bftheta$ in~\eqref{eq: lagrangian_coordinates_saddlepoint} we arrive at
\begin{equation}\label{eq : a_optimization}
        \inf_\bfa \medmath{\frac{1}{2} \int_0^T  }\bfa_s^\top\, \bfK_r^{-1}\bfa_s \,ds
        \,-\, \Ex_{\bfx\sim\rho_0} \phi_\bfa (0,\bfx),
\end{equation}
which is then the dual formulation of the optimization of~\eqref{eq: mean_field_control_problem_Kr} and reveals again the variational nature of the global interaction coefficients. 

Some remarks are in order. First, we have transformed the \emph{coupled} optimization problem~\eqref{eq: mfc_original_objective} into a saddle point problem where the agents are \emph{completely decoupled} at the expense of introducing a moderate number of dual variables (global interaction coefficients). 
That is, there is no longer a double expectation arising from~\eqref{eq: mean_field_interaction} to compute the interaction between every pair of agents. Importantly, for a fixed vector function $\bfa$, we can solve for $\bftheta$ in a completely parallel manner. 
In this case our original control $\bftheta$ is the primal variable and the interaction coefficients $\bfa_t$ are the dual variables.  
Second, the global interaction coefficients $\bfa_t$ can be \emph{reused}. More specifically, the primal-dual problem only needs to be solved \emph{once}. One may then solve another instance of~\eqref{eq: saddle_point_formulation}. In our experiments we show how \eqref{eq: saddle_point_formulation} can be solved for one instance using 100, 400, 800, and 1000 agents and can then reuse the corresponding $\bfa^\star_{100}, \bfa^\star_{400}, \bfa^\star_{800},$ and $\bfa^\star_{1000}$ for an experiment with a new instance of the MFC problem with 1000 agents (see Sec.~\ref{sec: recycling_a}).

\subsection{Discretization and Primal-Dual Algorithm}
To approximate the solution to~\eqref{eq: lagrangian_coordinates_saddlepoint}, we use the direct transcription approach~\cite{enright1992discrete} for simplicity. However, we note that any optimal control/trajectory generation algorithm~(e.g., \cite{zhao2022offline,ruthottomean2020,onkenneural}) can be used with our primal-dual framework when solving for $\bftheta$ in~\eqref{eq: lagrangian_coordinates_saddlepoint}.
Fix  a uniform time discretization
\begin{align}
    0=t_1 < t_2 < \cdots < t_{n} = T,
    \label{eq:tgrid}
\end{align}
with time-step $h=t_{k}-t_{k-1}$.
We define a discrete formulation where the time integrals in~\eqref{eq: lagrangian_coordinates_saddlepoint} are approximated using sums over the grid \eqref{eq:tgrid} and the expectations are approximated using sample averages over $N$ agents $\bfx_1,\ldots,\bfx_N$
drawn independently from $\rho_\circ$:
\begin{equation} 
    \begin{split}
    \label{eq: discretized_saddle_point_1}
    & \Lc(\bfa,\bftheta) = \medmath{\frac{h}{2}\sum_{k=1}^{n}} \,\bfa_{t_k}^\top\,\bfK_r^{-1}\,\bfa_{t_k}  \\
    &- \medmath{\frac{h}{N}\sum_{k=1}^n\sum_{\ell=1}^N}\,L(t_k, \bfz_{\bfx_\ell}(t_k), \bftheta_{\bfx_\ell}(t_k))  \\
    &\quad - \medmath{\frac{h}{N}\sum_{k=1}^n\sum_{\ell=1}^N}\,\bfa_{t_k}^\top\bfzeta(\bfz_{\bfx_\ell}(t_k)) - \medmath{\frac{1}{N}\sum_{\ell=1}^N}\,G(\bfz_{\bfx_\ell}(T)),
    \end{split}
\end{equation}
where $\bfz_{\bfx_\ell}(t_k)$ is the discretized state for the $\ell^{th}$ agent at time $t_k$ obtained  using, e.g., Euler's method to solve the ODE~\eqref{eq : dynamics}. From now on we will use $\bftheta$ and $\bfa$ to denote, respectively, the matrices with entries $\bftheta_{\bfx_\ell}(t_k)$ and $\bfa_{t_k}$, and 
$\bftheta_{\bfx} =(\bftheta_\bfx(t_1),\ldots,\bftheta_\bfx(t_n))$.

Importantly, we note again that $\mathcal{L}$ is \emph{separable} with respect to the agents; that is, we may write,
\begin{equation}
    \Lcurly(\bfa, \bftheta) = \medmath{\frac{1}{N}\sum_{\ell=1}^N}\,\, \Lc_\ell(\bfa, \bftheta_{\bfx_\ell}),
\end{equation}
where  
\begin{equation} 
    \begin{split}
    \label{eq: discretized_saddle_point}
    &\Lc_\ell(\bfa,\bftheta_\ell) = \medmath{\frac{h}{2}\sum_{k=1}^{n}} \,\bfa_{t_k}^\top\,\bfK_r^{-1}\,\bfa_{t_k}  \\
    &- \medmath{{h}\sum_{k=1}^n}\,L(t_k, \bfz_{\bfx_\ell}(t_k), \bftheta_{\bfx_\ell}(t_k))  \\
    &\quad - \medmath{{h}\sum_{k=1}^n }\,\bfa_{t_k}^\top\bfzeta(\bfz_{\bfx_\ell}(t_k)) - G(\bfz_{\bfx_\ell}(T)).
    \end{split}
\end{equation}

Following \cite{NurbekyanOne2018, Liucomputational2021}, we approximate the saddle point of~\eqref{eq: discretized_saddle_point} with the following primal-dual iterations
\begin{align}
    &\bftheta_{\bfx_\ell}^{k+1} = \argmax_{\bftheta} \,\Lc_\ell(\bfa^k, \bftheta_{\bfx_\ell}), \quad \ell = 1,\ldots, N \label{eq: primal_update}\\
    &\bfa^{k+1} = \argmin_\bfa\, \Lc(\bfa, \bftheta_\ell^{k+1} )+\frac{\|\bfa - \bfa^k\|^2}{2\gamma}, \label{eq: dual_update}
\end{align} 
where $\gamma >0$ is a chosen parameter and the optimization variables $(\bfa, \bftheta)$ are initialized randomly.
We note that the primal update~\eqref{eq: primal_update} can be done \emph{in parallel} for each agent and can be solved efficiently using, e.g., L-BFGS~\cite{liu1989limited}, while the dual update consists of minimizing a quadratic term with the analytic solution 
\begin{equation}
    \bfa^{k+1}_{t_i} = \left(\bfI + h_a \mathbf{K}_r^{-1} \right)^{-1} 
    \left( \bfa^{k}_{t_i}+h_a \medmath{\frac{1}{N} \sum_{l=1}^N}\, \bfzeta\left(\bfz_{\bfx_l}^{\bftheta_i^{k+1}}(t_i)\right) \right)
\end{equation}
for $i = 1, 2, \ldots, n$, where $h_a = h \gamma$.
In our experiments, $\mathbf{K}_r$ is the identity matrix because we use a kernel expansion of the form 
\begin{equation}
    K_r(\bfx,\bfy) = \bfzeta(\bfx)^\top \bfzeta(\bfy) ,
\end{equation}
where $\bfzeta \colon \mathbb{R}^d \to \mathbb{R}^r$ is a neural network (see Sec.~\ref{sec: kernel_approx})
    
\section{Related Works}

\subsection{Multiagent Control}
Traditional approaches to multiagent control do not operate in the mean-field framework.
The paper \cite{jadbadaiecoordination2003} introduces a theoretical framework for the discrete time model for autonomous agents presented in \cite{vickeknovel1995} that sparked  a variety of studies exploring multiagent control.
One can now find many surveys and tutorials on multiagent control; for example,  \cite{murrayrecent2007, gazicoordination2007, caooverciew2012, wileycooperative2007, wangdistributed2010, Bellinghammulti2003, lewiscooperative2014}.  
The work in \cite{bojinovmultiagent2002} provides a multiagent framework where agents are self-reconfigurable to aid in problems where adaptive structures are beneficial such as environmental stresses.
\cite{xiaofinite2009} provides an approach to multiagent control that splits the problem based on local and global interaction information in order to reduce the cost when the number of agents is large. \cite{dimarogonasevent2009} and \cite{ dimarogonastriggered2009} provide an event-triggered approach to multiagent control where the stability of a control scheme using agents with limited resources is explored.
The works in \cite{onkenneural,  verma2024neural, li2022neural, onken2021neural} propose methods for multi-agent optimal control that use a neural network to parameterize the value function. A survey on these approaches and their connections to optimal transport and generative modeling can be found in~\cite{ruthotto2024differential}. 
Several interesting applications of multiagent control include wind turbines \cite{johnsoncontrol2006}, vehicle coordination \cite{Kirchnerhamilton2020}, and reusable launch vehicles \cite{gongpseudospectral2007}. Other quadrotor simulation experiments can be found in \cite{Honigtrajectory2018, mengsympocnet2022}.  We refer to \cite{knorncollective2016}
for multiagent swarm planning using optimal control.

\subsection{Mean-Field Control}

Various studies have explored approaches to solving mulit-agent swarm problems using MFC. The early work \cite{yangmultiagent2008} provides a multiagent control approach for interacting agents related to MFC by incorporating statistical estimates of the global system properties.
\cite{zheng2021transporting} introduces an MFC approach that replaces local coordination strategies and is supplemented with stability analyses and velocity control experiments consisting of 1024 robot swarms.  \cite{tembine2013risk} introduces an MFC framework where the agents are risk-sensitive through the addition of exponential utility functions for independent agents. \cite{elamvazhuthi2018mean, elamvazhuthi2017mean} provide an MFC framework that uses continuous-time Markov chains for evolving agent states that emphasizes decentralized control laws for stabilizing agent distributions, and summarize numerical experiments for 80 and 1200 2D agents. 
The related control approach~\cite{elamvazhuthi2015optimal} uses stochastic strategies for swarms of resource-constrained robots by solving PDE-constrained optimization problems. 
\cite{ringh2023mean} proposes an MFC framework that allows for the use of different types of robot controllers in rescue scenarios when the dynamics for different robots are defined with different velocity-weighting strategies.  
The works~\cite{lee2021controlling, couillet2012electrical} and~\cite{emamiage2024} provide specific applications of MFC; namely, controlling the spatial velocity of epidemics, the behavior of a large number of electric and hybrid vehicle owners attempting to minimizing their operating cost, and age of information minimization using multi-agent robot swarms, respectively.  
For a comprehensive introduction/review of MFCs, we refer the reader to the survey \cite{elamvazhuthi2019mean} and references therein. 

Our proposed method is most closely related to \cite{Agrawal2022RandomFF}, which employs a similar kernel expansion approach but is used in the context of potential mean-field games (MFG).  In particular, \cite{Agrawal2022RandomFF} uses a random Fourier features approximation of the interaction term as opposed to a neural network approach. 
A similar method uses polynomial bases for low-dimensional MFGs~\cite{Nurbekyan2019Fourier}.  We refer to the survey \cite{wangsurvey2024} for MFG approaches that decouple interaction and agent control.
Our work is also motivated by the works in~\cite{ruthottomean2020, lin2021alternating} and \cite{wang2021global}; the former uses deep learning to solve MFC and MFG problems.  The latter combines policy gradient estimation with the MFC framework and provides a linear convergence result. Hwever, neither work considers non-local interactions. 

\section{Numerical Experiments}
\label{sec: numerical_experiments}
In order to demonstrate the effectiveness of our primal-dual framework, we conduct experiments with swarms consisting of 1000 and 5000 agents. The primal-dual formulation proves to be valuable in addressing challenges associated with large-scale swarm control where the dimensionality of the problem arising from the number of agents poses significant computational hurdles when computing the non-local interaction between agents. 
The generated trajectories showcase coordinated motion of interacting agents using two types of dynamics, double integrator and quadcopter~\cite{carrillo2012quad}.  
These experiments emphasize the algorithm's ability to numerically solve optimal control of swarms.
For all experiments, we choose the following stopping criteria:
$\| \nabla_{\bftheta} \mathcal{L}(\bftheta^k) \| \leq \epsilon_{tol}$ (optimality criterion for primal variable), $\big\| \bfa^{k}_{t_i} - \medmath{\frac{1}{N} \sum_{l=1}^N} \, \bfzeta(\bfz_{\bfx_l}^{\bftheta_i^{k+1}}(t_i)) \big\| < \epsilon_{tol}$ for all $t_i$ (optimality criterion for dual variable), and $\| \nabla_\theta \mathcal{J}_r(\theta^k) \| < \epsilon_{tol}$ in order to ensure the MFC problem~\eqref{eq: mean_field_control_problem_Kr} is solved. Here, we choose $\epsilon_{tol} = 0.5$.

\subsection{Kernel Approximation}
\label{sec: kernel_approx}

In our experiments, we use exponential interaction kernels
\begin{equation}\label{eq:exp_kernel}
    K(\bfx,\bfy)=\alpha_1 e^{-{\|\bfx-\bfy\|^2}/{2}},\quad \bfx,\bfy \in \mathbb{R}^3.
\end{equation}
We build the approximation in~\eqref{eq: kernel_expansion_framework} in two steps. First, we find $\Tilde{\bfzeta}=(\Tilde{\zeta}_1,\Tilde{\zeta}_2,\cdots,\Tilde{\zeta}_r)$ such that
\begin{equation}\label{eq:net_training}
    \Tilde{K}(\bfx,\bfy)=e^{-{\|x-y\|^2}/{2}}\approx \Tilde{\bfzeta}(\bfx)^\top \Tilde{\bfzeta}(\bfy).
\end{equation}
Next, we set
\begin{equation*}
    \bfzeta(\bfx)=\sqrt{\alpha_1} \,\Tilde{\bfzeta}(\bfx).
\end{equation*}
Hence, we achieve~\eqref{eq: kernel_expansion_framework} with $\bfK_r$ being the identity matrix.

\begin{figure}[t] 
    \centering
	\begin{tabular}{cc}
        $K$ & $K_r$
        \\
        \includegraphics[width=0.2\textwidth]{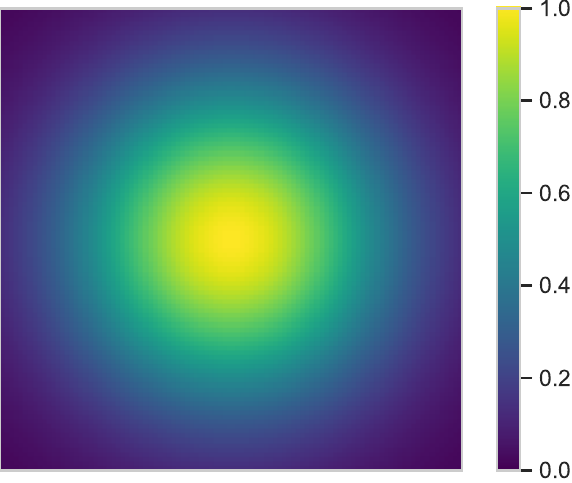}&
        \includegraphics[width=0.2\textwidth]{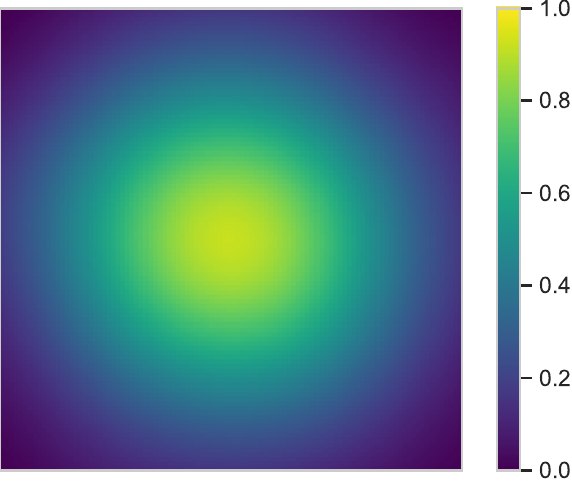}
	\end{tabular} 
\caption{\textbf{Kernel Approximation: } Gaussian kernel (left) and neural network kernel approximation using $r=50$ basis functions (right).}
  \label{app: approximatekernel}
\end{figure}

{For~\eqref{eq:net_training}, we train a single-layer feed-forward neural network $\Tilde{\bfzeta}$ with one 100-dimensional hidden layer. The network outputs $r=50$ basis functions $\Tilde{\zeta}_i$ for $i=1,\ldots, r$ and has 5450 total network parameters. The network is trained on $1e5$  samples drawn uniformly from the cube $[-3,3]^3$ 
for $5e4$ iterations using the Adam optimizer~\cite{kingmaadam2017} with a learning rate of $1e-3$ and a scheduler that sales the learning rate by $1e-1$ every $1e4$ iterations.  
Training occurs by minimizing the sample mean squared-error between the true kernel values and our approximate kernel values generated by the network. 
In order to obtain higher accuracy, we also penalize the norm of the difference of the true gradient and the approximate gradient of the kernel. 
In this experiment, we obtain a validation loss of 7.70e-03, which we found to be accurate enough for our MFC experiments. 
As an illustration, we show the true Gaussian kernel and its neural  network approximation in Figure~\ref{app: approximatekernel}.}

\subsection{Experimental Setup}
\label{sec: accel_obstacle_trajectory_experiment}
\subsubsection{Double Integrator}
We approximate the optimal controller of a swarm of agents with dynamics given by
\begin{equation}
\medmath{
\Dot{\mathbf{z}} =
\begin{cases}
\Dot{x} = v_x\\
\Dot{y} = v_y\\
\Dot{z} = v_z\\
\dot{v_x} = a_x\\
\dot{v_y} = a_y\\
\dot{v_z} = a_z\\
\end{cases}
}
\end{equation}
where the controls are denoted by the variable $\bftheta = (\Dot{v_x}, \Dot{v_y}, \Dot{v_z})\in \R^3$ and $\Dot{v}$ are the accelerations that correspond to each spatial dimension.
%
\begin{figure*}[t]
    \centering
    \begin{tabular}{cccc}
         \multicolumn{2}{c}{\textbf{1000 Agents: Double Integrator}}&  \multicolumn{2}{c}{\textbf{1000 Agents: Quadrotors}} \\
        \includegraphics[width=0.2\textwidth]{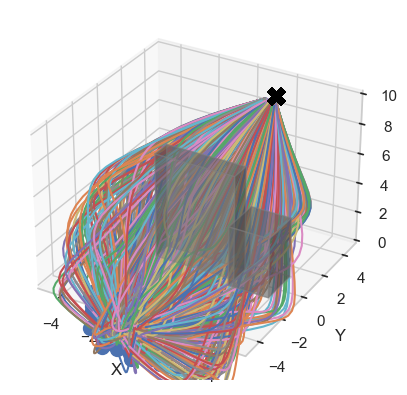}
        & 
        \includegraphics[width=0.24\textwidth]{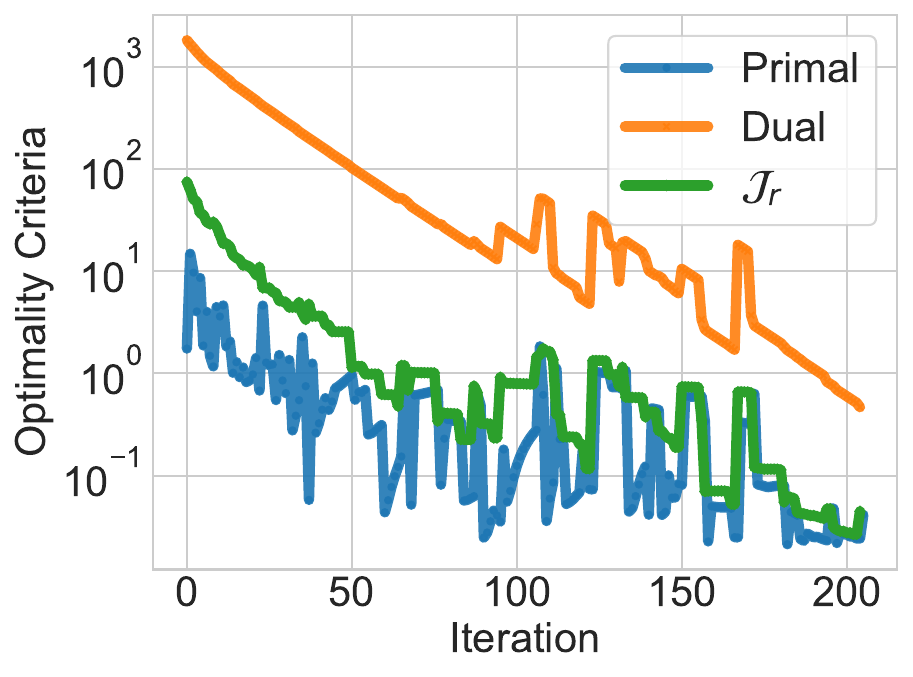}
        &
        \includegraphics[width=0.2\textwidth]{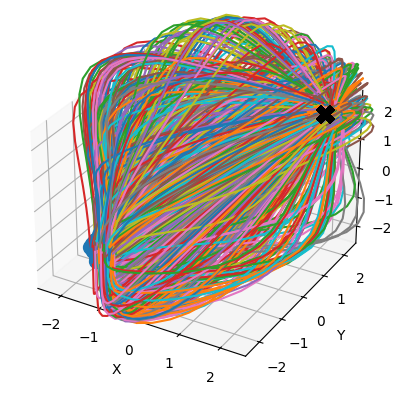}
        & 
        \includegraphics[width=0.24\textwidth]{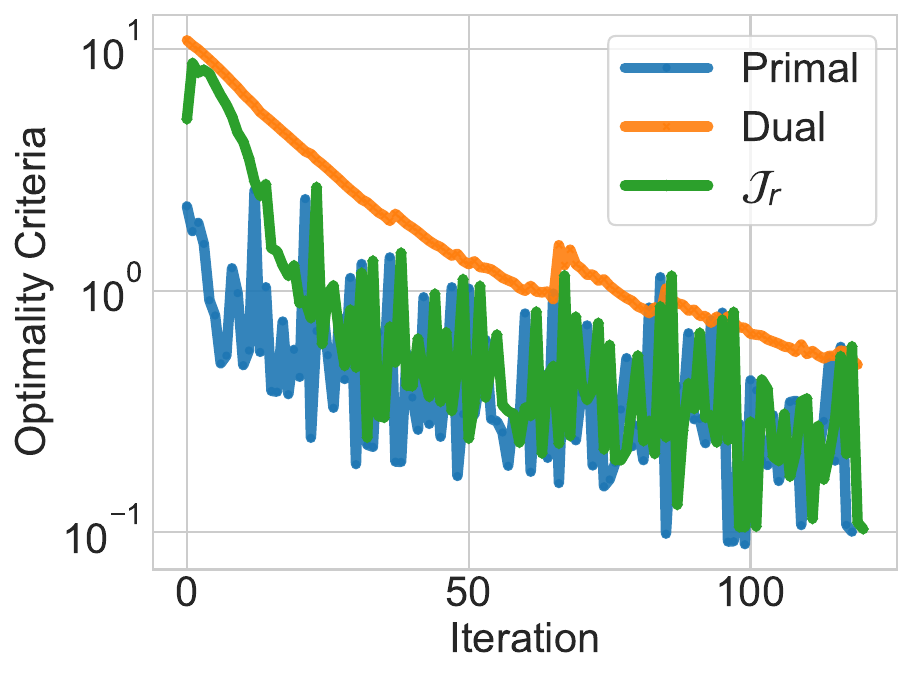}
        \\
    \end{tabular}
    \caption{Swarm Control Experiments. \textbf{Left panel:} Trajectories and gradient norms for double integrator dynamics with 1000 agents. \textbf{Right panel:} Trajectories and gradient norms for quadrotor dynamics with 1000 agents. On the left of each panel, trajectories are shown where each agent must reach the target (black X) and avoid colliding with one another. On the right of each panel, gradient norm history of primal-dual is shown: The primal and dual optimality criterion are plotted in blue and orange, respectively. The gradient norm of the original MFC problem~\eqref{eq: mean_field_control_problem_Kr} with the approximate kernel is shown in blue.}
    \label{fig:combined_1000_doubletinegrator_quadrotor_experiments}
\end{figure*}
{Agents  are  penalized if they get too close to each other using the interaction kernel~\eqref{eq:exp_kernel} and the penalty parameter $\alpha_1 = 2e5$ to weigh the interaction cost.}

The agents must also avoid colliding into two obstacles.
The obstacles are modeled as a function $Q:\R^3\to \R$ defined as in \cite{onkenneural}; namely,
\[
Q(\bfx) = Q_A(\bfx) \mathbf{1}_{A}(\mathbf{x})+Q_B(\bfx) \mathbf{1}_{B}(\mathbf{x}),
\]
where  
$A = \medmath{[-2,2] \times [-0.5, 0.5] \times [0, 7]}$ and $B = \medmath{[2,4]\times [-1,1] \times [0,4]}$,
$\mathbf{1}_{S}$ is the indicator function of a set $S$,
and $Q_A$ and $Q_B$ are, respectively, the density functions of the Gaussian distributions $N(\bfmu_1,\bfSigma_1)$ and $N(\bfmu_2,\bfSigma_2)$ with 
$\bfmu_1 = (0,0,2)$, $\bfmu_2 = (2.5, 0, 2)$, $\bfSigma_1 = \mathrm{Diag}\{9, 3, 9\}$, and
$\bfSigma_2 = \mathrm{Diag}\{9,3,3\}$. To avoid the obstacles, $Q$ is included in the Lagrangian as
\begin{equation}
    L(t,\bfz_\bfx(t),\bftheta) = \|\bftheta_\bfx(t)\|^2 + \alpha_2 \,Q(\bfz_\bfx(t))
    \label{eq: lagrangian_double}
\end{equation}
with $\alpha_2 = 1e7$.
Finally, we use the terminal cost
 \begin{equation}
    G(\bfz_\bfx(T)) = \medmath{\frac{\alpha_3}{2}}\,\|\bfz_\bfx(T) - \bfz_\mathrm{target}\|^2
    \label{eq: terminal_cost}
\end{equation}
with $\alpha_3 = 1e4$.
The agents are initialized randomly as
$\bfz_{\bfx_l}(0) \sim {N}(\bfx_{\text{mean}}, 0.8 \,\boldsymbol{I})$, $l=1,2,\ldots,N$, and $\bfx_{\text{mean}} =(0, -0.5, 0, \ldots,0)\in \R^{6}$.
The target for the agents is chosen to be the point $(0,0,7,0,\ldots, 0) \in \R^{6}$.

To solve the primal problem, we use 250 L-BFGS~\cite{liu1989limited} iterations for fast initial improvement followed by 750 iterations of gradient descent to ensure convergence of the gradient norm. 

\subsubsection{Quadrotor}
We approximate the optimal controller of a swarm of quadrotors whose dynamics are given by 
\begin{equation}
\medmath{
\Dot{\mathbf{z}} = 
    \begin{cases}
        \Dot{x} = v_x\\
        \Dot{y} = v_y\\
        \Dot{z} = v_z\\
        \Dot{\psi} = v_{\psi}\\
        \Dot{\theta} = v_{\theta}\\
        \Dot{\phi} = v_{\phi}\\
        \Dot{v}_x = \frac{u}{m}\left(\sin(\phi)\sin(\psi) + \cos(\phi)\cos(\psi)\sin(\theta)\right)\\
        \Dot{v}_y = \frac{u}{m}\left(-\cos(\psi)\sin(\phi) +\cos(\phi)\sin(\theta)\sin(\psi)\right)\\
        \Dot{v}_z = \frac{u}{m}\cos(\theta)\cos(\phi) - g\\
        \Dot{v}_{\psi} = \Tilde{\tau}_{\psi}\\
        \Dot{v}_{\theta} = \Tilde{\tau}_{\theta}\\
        \Dot{v}_{\phi} = \Tilde{\tau}_{\phi}
    \end{cases}
    }
\end{equation}
where the state\\
$\bfz = (x,y,z,\psi,\theta,\phi, v_x, v_y, v_z, v_{\psi}, v_{\theta}, v_{\phi})\in \R^{12}$ is comprised of the spatial positions $(x,y,z)$, the angular orientation of the quadrotor, $(\psi, \theta, \phi)$, and $v$ defines the velocities associated with the spatial positions and angular orientations.  The controls are defined by the variable $\bftheta = (u, \Tilde{\tau}_{\psi}, \Tilde{\tau}_{\theta}, \Tilde{\tau}_{\phi}) \in \R^4$ where $u$ is the primary propulsion force oriented downward from the underside of the quadrotor and $(\Tilde{\tau}_{\psi}, \Tilde{\tau}_{\theta},\Tilde{\tau}_{\phi})$ define the torques corresponding to $(\psi, \theta, \phi)$. 
In this setup, agents must avoid collision with one another. The penalty parameters are the same as those described in the double-integrator experiments and are chosen as $\alpha_1 = 5e4$, $\alpha_2 = 0$, and $\alpha_3 = 2e3 $. In each primal-dual iteration, we use 20 L-BFGS~\cite{liu1989limited} iterations for the primal problem~\eqref{eq: primal_update} per dual update~\eqref{eq: dual_update}. 

 \begin{figure*}[t] 
    \centering
	\begin{tabular}{cccc}
        $\bfa^{\star}_{100}$ & $\bfa^{\star}_{400}$ & $\bfa^{\star}_{800}$ & $\bfa^{\star}_{1000}$\\
        \includegraphics[width=0.24\textwidth]{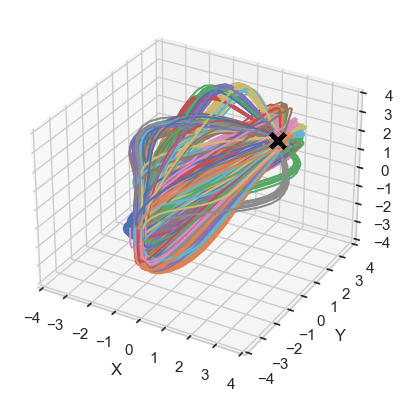}&
        \includegraphics[width=0.24\textwidth]{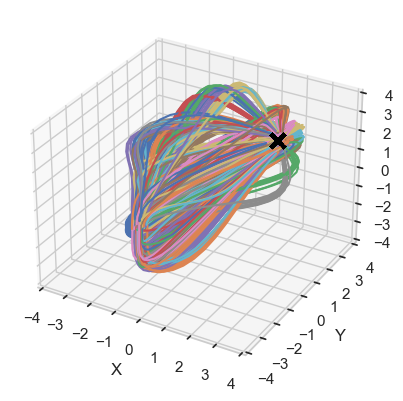}&\includegraphics[width=0.24\textwidth]{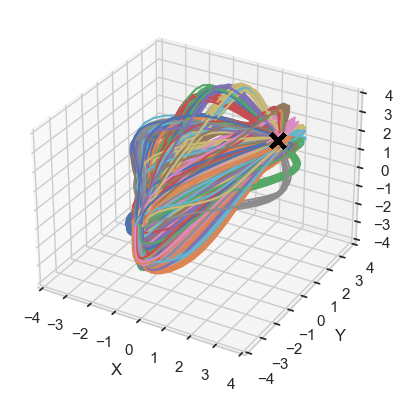}&
        \includegraphics[width=0.24\textwidth]{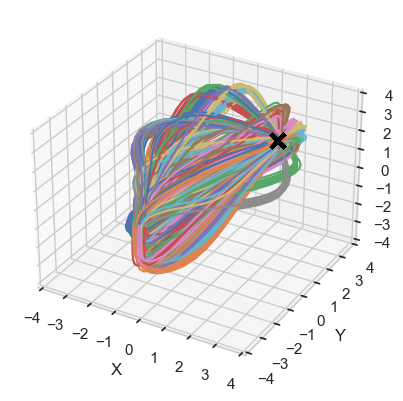}\\

        $\frac{\|\bfz_{\bfa^\star_{100}}- \bfz_{\bfa^\star_{1000}}\|}{\| \bfz_{\bfa^\star_{1000}}\|} $ & $\frac{\|\bfz_{\bfa^\star_{400}}- \bfz_{\bfa^\star_{1000}}\|}{\| \bfz_{\bfa^\star_{1000}}\|} $ & $\frac{\|\bfz_{\bfa^\star_{800}}- \bfz_{\bfa^\star_{1000}}\|}{\| \bfz_{\bfa^\star_{1000}}\|} $ &$\frac{\|\bfz_{\bfa^\star_{1000}}- \bfz_{\bfa^\star_{1000}}\|}{\| \bfz_{\bfa^\star_{1000}}\|} $\\
        \\[0.1mm]
        $0.0302$ & $0.0301$ & $0.0304$ & $0.000$ \\
	\end{tabular} 
	\centering 
	\\[2mm]
\caption{\textbf{Quadrotor Controller: } Illustration of trajectories using different mean-field interaction coefficients. 
Here, $a^\star_{100}, a^\star_{400}, a^\star_{800}$, and $a^\star_{1000}$ are interaction coefficients by solving the primal-dual for {100, 400, 800, and 1000} agents, respectively. Next, to make a proper comparison among the different coefficients, we sample a fixed set of 1000 initial conditions from $\rho_0$ and re-solve the primal problem using $a^\star_{100}, a^\star_{400}, a^\star_{800}, a^\star_{1000}$.
}
  \label{fig:a_ceoff_trajectory_experiment}
\end{figure*}

\subsection{Results}
For both dynamics, we plot trajectories and optimality criteria using 1000 agents in Figure~\ref{fig:combined_1000_doubletinegrator_quadrotor_experiments}.
For optimality criteria, we plot the gradient norm history of the objective function~\eqref{eq: mean_field_control_problem_Kr} (in green) along with the primal and dual optimality criteria (in blue and orange). 

Our results show that our proposed primal-dual approach solves the proposed problem~\eqref{eq: mean_field_control_problem_Kr} \emph{in a decoupled manner} as seen in the iterates~\eqref{eq: primal_update}.
The 5000-quadrotor experiment shown in Figure~\ref{fig: 5k_quadcopter_trajectories} in the introduction has the same experimental setup. These experiments were run on an Apple MacBook Pro 8-core 3200 MHz M1 laptop running macOS Monterey 12.0.1 with a total of 16GB of RAM.

\subsection{Time Comparison: Coupled vs. Primal-Dual}
To illustrate why the primal-dual approach is preferred, we compare the runtimes between solving~\eqref{eq: original_MFC_problem} directly (we refer to this as the coupled approach) and the proposed primal-dual approach using the quadrotor dynamics with 20 and 50 agents.
\begin{figure}[ht]
    \centering
    \begin{tabular}{cc}
    \hspace{10mm}\textbf{20 Agents} & \hspace{10mm}\textbf{50 Agents}\\
    \includegraphics[width=.23\textwidth]{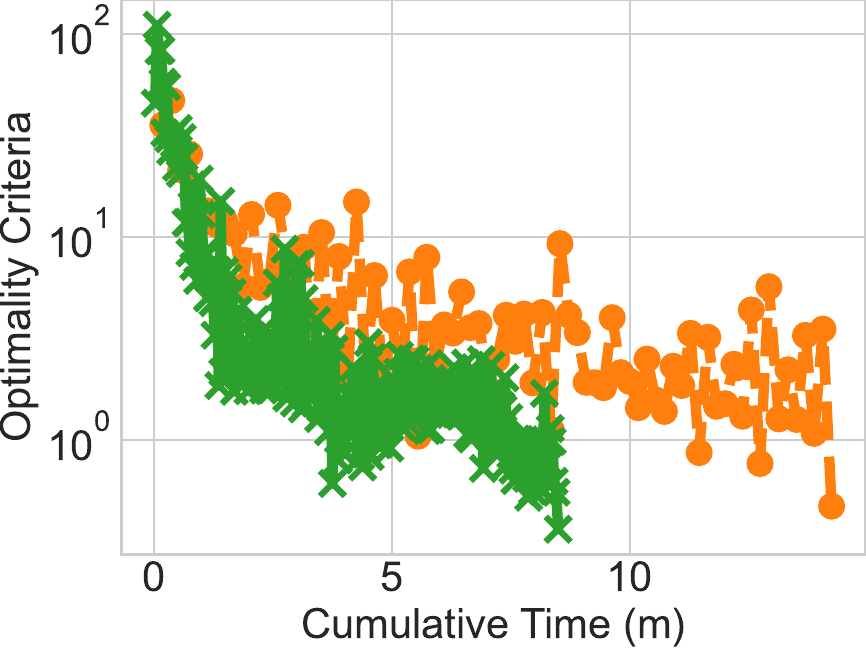} 
    & \includegraphics[width=.23\textwidth]{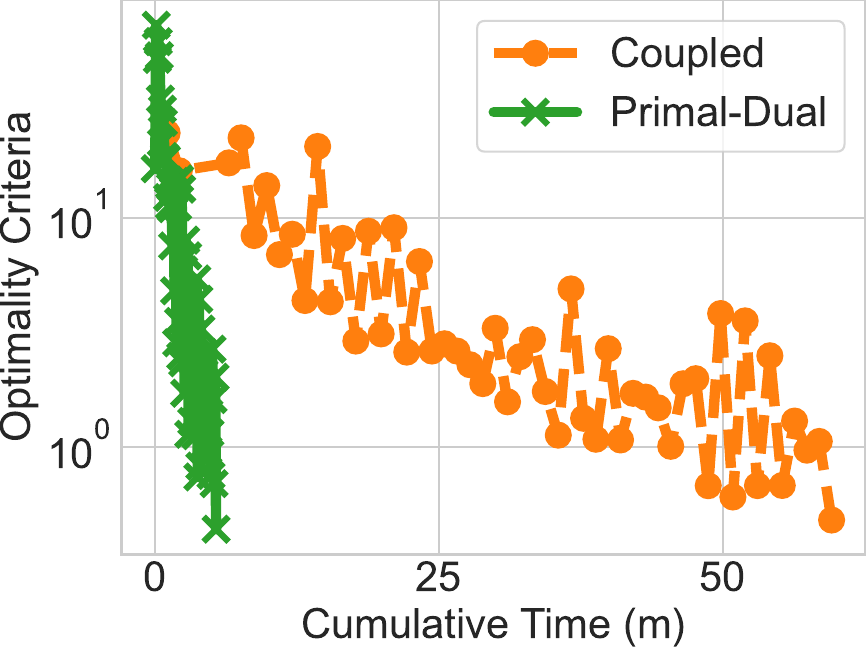}\\
    \end{tabular}
    \caption{Gradient norm vs. cumulative runtime (minutes) for the objective function $J_r$ in~\eqref{eq: mean_field_control_problem_Kr} for 20 and 50 agents. Each algorithm is run three times and the fastest run is chosen. Optimization stops when the gradient norm reaches $5e-1$. Orange shows the coupled approach and green shows our primal-dual approach.}
    \label{fig:quadrotor_gradnorm_vs_totaltime_dloop_and_broadcasting}
\end{figure}
Figure~\ref{fig:quadrotor_gradnorm_vs_totaltime_dloop_and_broadcasting} shows the advantages of using the primal-dual approach, especially as the number of agents increases. 
The experiments were run on an Apple MacBook Pro 8-core 3200 MHz M1 laptop running macOS Monterey 12.0.1 with a total of 16GB of RAM. As expected, the primal-dual approach is faster in terms of cumulative runtime due to the sequential nature of the coupled approach. It is worth noting that our simple implementation did not fully utilize the parallel nature of our primal-dual framework and further time gains could be achieved with a more advanced parallel implementation (e.g., using a star network topology infrastructure~\cite{zhang2014asynchronous}).

\subsection{Reusing the Global Mean-Field Interaction}
\label{sec: recycling_a}
An important benefit of our primal-dual approach is that the global interaction coefficients $\bfa$  can be computed, saved, and then reused on future instances of the MFC problem~\eqref{eq: original_MFC_problem}. This is because once $\bfa$ is obtained, only the primal problem requires optimization \emph{regardless of the number of agents}.
Consequently, $\bfa$ can be determined using \emph{fewer} agents than the number of trajectories that are required as long as the interaction is adequately approximated.
This is particularly beneficial for situations where, for
example, some quadrotors might lose power because re-computing the interaction is \emph{not required}. We empirically demonstrate this in Figure~\ref{fig:a_ceoff_trajectory_experiment} as follows. First, we compute the optimal coefficients $a_{100}^\star$ using 100 agents, $a_{400}^\star$ using 400 agents, 
$a_{800}^\star$ using 800 agents, and 
$a_{1000}^\star$ using 1000 agents. 
Second, we sample 1000 fixed initial conditions $\bfx \sim \rho_0$. Finally, we re-solve the \emph{primal problem only} with these new initial points using $\bfa_{100}, \bfa_{400}, \bfa_{800},$ and $\bfa_{1000}$ until the primal optimality condition is satisfied.
{ Figure~\ref{fig:a_ceoff_trajectory_experiment} shows that, as long as there are enough samples, trajectories generated using $a^\star_{100}, a^\star_{400}, a^\star_{800},$ and $a_{1000}^\star$ are qualitatively and quantitatively similar. Finally, as previously stated, Figure~\ref{fig:a_ceoff_trajectory_experiment} indicates that one may solve for $a^\star$ using \emph{fewer} agents and use this pre-computed $a^\star$ to solve new instances of the MFC problem, with a potentially larger number of agents, and in a completely \emph{parallel manner}.}


\begin{figure}[t]
    \centering
    \includegraphics[width=0.3\textwidth]{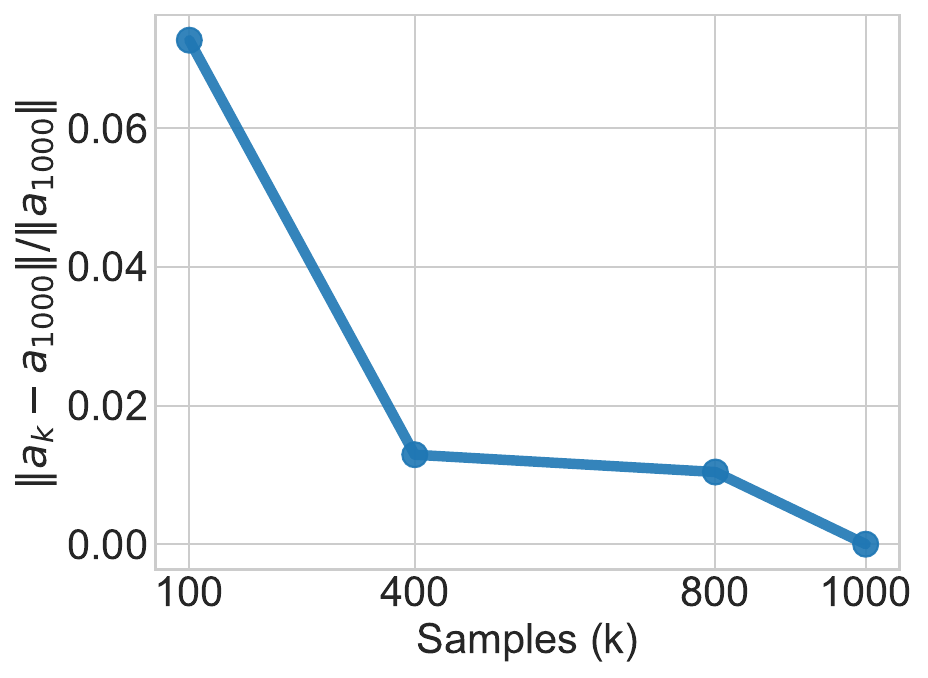} \\
    \caption{Number of agents used to train the mean-field interaction coefficients $\bfa$ vs. relative norm difference $\|\bfa_k - \bfa_{1000} \| / \|\bfa_{1000}\|$ for $k=100$, $400$, $800$, and $1000$ agents for quadrotor dynamics.}      
    \label{fig:quadrotor_acoeff_normdifferences}
\end{figure}


\section{Conclusions}
We propose a primal-dual framework based on kernel expansions for solving mean-field control (MFC) problems. Our approach decouples agent interactions by introducing global interaction coefficients based on kernel expansions. These coefficients allow us to reformulate the MFC problem as a saddle point problem where the primal problem can be solved in parallel. 
Moreover, new instances of the MFC problem can be solved at reduced costs by reutilizing the global interaction coefficients. 
Our experiments show that the primal-dual framework is effective at solving MFC problems, including an optimal swarm control problem of 5000 quadrotors.
As there are evident connections between optimal control and optimal transport~\cite{ruthottomean2020, lin2021alternating}, future works include using similar techniques in the context of generative modeling~\cite{vidaltaming2023, OnkenOTFlowFA2020, finlay2020train, ruthotto2021introduction, wang2023efficient} and inverse MFCs~\cite{chow2022numerical, liu2023inverse, ding2022mean}. More elaborate primal-dual solvers such as alternating direction method of multipliers~\cite{boyd2011distributed, fung2019uncertainty, ye2022adaptive, fung2020admm} will also be explored.

\bibliographystyle{abbrv}      
\bibliography{references.bib}

\section{}

\section*{Acknowledgments}
Samy Wu Fung and Alexander Vidal were  partially funded by National Science Foundation award DMS-2110745. Stanley Osher was partially funded by Air Force Office of Scientific Research (AFOSR) MURI FA9550-18-502 and Office of Naval Research (ONR) N00014-20-1-2787.

\ifCLASSOPTIONcaptionsoff
  \newpage
\fi

\end{document}